\documentclass[]{interact}

\usepackage{epstopdf}% To incorporate .eps illustrations using PDFLaTeX, etc.
\usepackage[caption=false]{subfig}% Support for small, `sub' figures and tables

\usepackage[numbers,sort&compress]{natbib}% Citation support using natbib.sty
\bibpunct[, ]{[}{]}{,}{n}{,}{,}% Citation support using natbib.sty
\renewcommand\bibfont{\fontsize{10}{12}\selectfont}% Bibliography support using natbib.sty

\theoremstyle{plain}% Theorem-like structures provided by amsthm.sty
\newtheorem{theorem}{Theorem}[section]
\newtheorem{lemma}[theorem]{Lemma}

\theoremstyle{definition}

\theoremstyle{remark}
\newtheorem{remark}{Remark}[section]

\usepackage{hyperref}
\usepackage{graphicx}
\usepackage{dsfont}
\usepackage{mathtools}
\usepackage{wrapfig}
\usepackage{amssymb, amsmath, latexsym}
\usepackage{nicefrac}
\usepackage{hhline}
\usepackage{multirow}
\usepackage{pifont}% http://ctan.org/pkg/pifont
\usepackage[colorinlistoftodos,bordercolor=orange,backgroundcolor=orange!20,linecolor=orange,textsize=scriptsize]{todonotes}
\usepackage{algorithm}\usepackage{algpseudocode}
\usepackage[symbol]{footmisc}

\usepackage{tikz}
\usetikzlibrary{matrix, positioning, fit}

\newcommand{\argmin}{\mathop{\arg\!\min}}

\newcommand{\Var}{\mathrm{Var}}

\newcommand*\Laplace{\mathop{}\!\mathbin\bigtriangleup}

\begin{document}

% \articletype{ARTICLE TEMPLATE}% Specify the article type or omit as appropriate

\title{Gradient directions and relative inexactness in optimization and machine learning\thanks{}}
%

% \titlerunning{Stopping rules for accelerated gradient methods with additive noise in gradient}
% If the paper title is too long for the running head, you can set
% an abbreviated paper title here
%

\author{
\name{
Artem Vasin \textsuperscript{a}
}
\affil{
\textsuperscript{a}Moscow Institute of Physics and Technology, Dolgoprudny, Russia;
}
}

\date{Received: date / Accepted: date}
\maketitle              % typeset the header of the contribution
\begin{abstract}
In this paper, we investigate the influence of noise giving an estimate of the gradient having a acute angle with the original. Noise amplitude has a relative model. The work offers both theoretical calculations and theorems, as well as experimental results. Classic machine learning problems were chosen as experiments - linear and logistic regression, computer vision and natural language processing.
\end{abstract}

% \clearpage

\section{Introduction}\label{Introduction}
We consider global optimization problem:
\begin{equation}
\label{optim}
    \min\limits_{x \in \mathbb{R}^n} f(x).
\end{equation}
We define $f^*$ - as minimum value of $f$ or solution for problem~\ref{optim} and $x^*$: $f(x^*) = f^*$, also for iterative methods with starting point $x^0$ we can define $R =  \|x^0 - x^* \|_2$.
We assume that the objective $f$ is $L$-smooth i.e., for all  $x,y\in \mathbb{R}^n$:
\begin{equation} \label{smooth cond}
    \begin{gathered}
        \|\nabla f(y) - \nabla f(x)\|_2 \leqslant L\|y - x\|_2, \\
        \text{Or equivalent and more usefull:}, \\
        f(y) \leqslant f(x) + \langle \nabla{f}(x), y - x \rangle + \frac{L}{2} \|x - y \|_2^2
    \end{gathered}
\end{equation}
Also we consider stochastic optimization:
\begin{equation} \label{stoh opt}
    \begin{gathered}
        f(x) \to \min\limits_{x \in \mathbb{R}^n}, \\
        \mathbb{E} \left[ \widetilde\nabla{f}(x) | x \right] = \nabla f(x).
    \end{gathered}
\end{equation}
The classical formulation of the machine learning problem is the sum-structured problem:
\begin{equation} \label{sum opt}
    f(\theta) = \frac{1}{m} \sum\limits_{i = 1}^m f_i(\theta) \to \min\limits_{\theta \in \mathbb{R}^n}.
\end{equation}
Here $\theta$ - model parameters and $f_i$ - loss function on sample $i$. Stochastic optimization in this problem is defined as follows:
\begin{equation*}
    \begin{gathered}
        I \sim \mathcal{U}\lbrace 1 \dots m \rbrace \times \dots \times \mathcal{U}\lbrace 1 \dots m \rbrace, \\
        \widetilde{f}(\theta) = \frac{1}{B} \sum\limits_{i \in I} f_i(\theta).
    \end{gathered}
\end{equation*}
Where $B$ - batch size. We introduce two model errors in gradient:
\begin{equation}\label{abs_inexact}
    \|\tilde{\nabla} f(x) - \nabla f(x)\|_2 \leqslant \delta, \text{ (absolute error) or }
\end{equation}
\begin{equation}\label{relative_inexact}
    \|\tilde{\nabla} f(x) - \nabla f(x)\|_2 \leqslant \alpha\|\nabla f(x)\|_2,  \text{ (relative error)}.
\end{equation}
For~\ref{relative_inexact} model we will use the following condition:
\begin{equation}\label{correl}
   (\forall x \in \mathbb{R}^n ): \; \langle\widetilde{\nabla} f(x), \nabla f(x) \rangle \geqslant \gamma \|\widetilde{\nabla}f(x) \|_2 \|\nabla{f}(x) \|_2, \quad \gamma \in (0, 1 ]
\end{equation}
We can interpret the condition~\ref{correl} as lower bound for cosine angle between gradient and it is estimation. Follow~\cite{vaswani2019fast, beznosikov2024first} we can define $(\nu, \rho)$ -- noise growth condition:
\begin{equation} \label{growth condition}
    (\forall x \in \mathbb{R}^n ): \; \nu \|\nabla f(x) \|_2 \leqslant \| \widetilde{\nabla} f(x) \|_2 \leqslant \rho \| \nabla{f}(x) \|_2
\end{equation}
We can note, that relative model~\ref{relative_inexact} with $\alpha \in [0; 1)$ implies~\ref{correl} and growth model~\ref{growth condition} with:
\begin{equation} \label{alpha to growth}
    \gamma = \sqrt{1 - \alpha^2}, \\
    \nu = 1 - \alpha, \\
    \rho = 1 + \alpha.
\end{equation}
We propose studies of the convergence of first-order methods with conditions~\ref{smooth cond}, \ref{correl}, \ref{growth condition}. We will prove, that classic gradient procedure~\ref{alg gd} will preserves the order of convergence up to constants:
\begin{equation*}
     f(x^N) - f^* = O \left( \frac{\rho^2}{\nu^2 \gamma^2} \frac{LR^2}{N} \right)
\end{equation*}
In Sections~\ref{section stm}, \ref{CG} provided motivation and relevant to \cite{kornilov2023intermediate}, \cite{gannot2021frequency} results, associated with $\sqrt{\frac{\mu}{L}}$, where $\mu$ - constant of strong convexity~\ref{strong mu}.

Paper contains a sufficient number of experiments with modern deep learning models with coefficients $\alpha, \gamma$ estimation. Sufficient conditions for the dataset are also given that guarantee the conditions~\ref{relative_inexact}, \ref{correl}.

\section{Ideas behind the results}\label{Ideas} 
Most papers consider absolute model~\ref{abs_inexact}, but what $\delta$ should we choose for theoretical estimation of convergence. For example Algorithm~\ref{alg gd} has convergence (for convex function with~\ref{smooth cond}):
\begin{equation*}
    f(x^N) - f^* = O \left( \frac{LR^2}{N} + \delta \right).
\end{equation*}
If estimation for $\delta$ is large the theoretical convergence will be uninformative, but if we plot convergence plot we will see decreasing graph. As an example we can take dataset CIFAR-10~\cite{krizhevsky2009learning} for classification problem. Dataset consist of 50k training samples and 10k test samples of 32x32 images with 10 classes. We will use ResNet-18~\cite{he2016deep} as classification model and PyTorch framework~\cite{Paszke_PyTorch_An_Imperative_2019}, because it provides batching. We will estimate $\alpha, \delta$ and $\gamma$ coefficient on each iteration by transforming epochs to iterations. Iterate over all batches (dataloader in PyTorch) we can summarize gradients per batches to gradient for whole train dataset, then choosing single batch we can evaluate required values.
\begin{figure}[H] \label{motiv graph}
	\center{\includegraphics[width=1\linewidth]{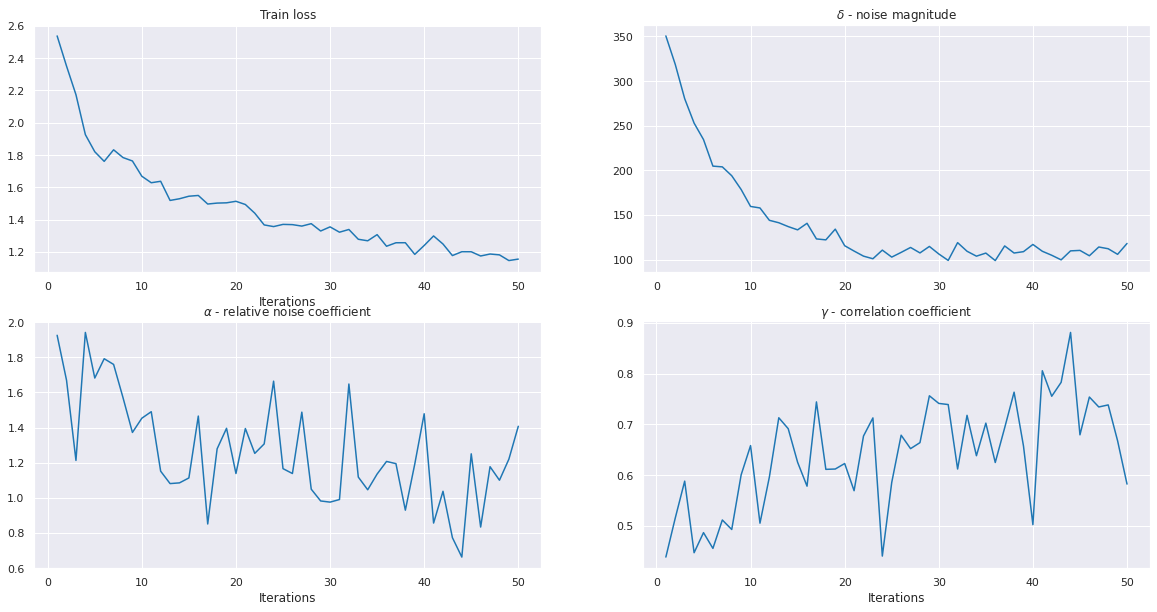}}
	\caption{Convergence and coefficient evaluation per iterations for ResNet-18 CIFAR-10.}
\end{figure}
We can see, that $\gamma$ - 0.45, gives intuition to explore convergence with such conditions. More experiments provided at Section~\ref{section exp}.
We should note, that for stochastic optimization~\ref{stoh opt} convergence can be much better, using only $\delta_* = \mathbb{E}\|\widetilde{\nabla}f(x^*) \|_2$, more details at~\cite{woodworth2020minibatch}.

\section{Motivation for relative noise}\label{motivation}
Most works consider an absolute noise model~\ref{abs_inexact}, for example ~\cite{devolder2014first, beznosikov2020gradient, dvurechensky2016stochastic}. However there are many papers consider relative model~\ref{relative_inexact} - \cite{gannot2021frequency, vaswani2019fast, kornilov2023intermediate, vasin2023accelerated}. These noise models can be approached from the point of view of stochastic differential equations, since the origin of the standard gradient procedure can be perceived as a discretization corresponding to an autonomous system.
For models~\ref{abs_inexact}, \ref{relative_inexact}, consider the following autonomous SDEs, respectively.
\begin{equation}\label{abs sde}
    d X_t = - \nabla{f}(X_t) dt + \delta \cdot d W_t,
\end{equation}
\begin{equation}\label{rel sde}
    d X_t = - \nabla{f}(X_t) dt + \alpha \|\nabla{f}(X_t) \|_2 \cdot d W_t.
\end{equation}
For SDE:
\begin{equation*}
    d X_t = b(X_t)dt + \Sigma(X_t) dW_t,
\end{equation*}
we can define infinitesimal operator:
\begin{equation} \label{infinitesimal}
    \mathcal{A} \phi = L_{b}\phi + \frac{1}{2} tr \left(\Sigma \Sigma^T \frac{\partial^2 \phi}{\partial x^2} \right), \quad L_b \phi = \langle b(x), \nabla \phi(x) \rangle.
\end{equation}
One can show:
\begin{equation*}
    \mathcal{A}^* \phi = -div(b \cdot \phi) + \frac{1}{2} \sum\limits_{i, j} \partial_{x_i, x_j} D_{i, j} \phi, \quad D(x) = \Sigma(x) \Sigma^T(x).
\end{equation*}
For the distribution solution of the stochastic differential equation, the Fokker-Planck equation is valid:
\begin{equation}\label{FP}
    \partial_t P = \mathcal{A}^* P.
\end{equation}
According to~\cite{xu2018global} stationary solution of Fokker-Planck equation~\ref{FP} for~\ref{abs sde}:
\begin{equation*}
    p^*(x) \propto \exp\left({-\frac{2 f(x)}{\delta^2}}\right).
\end{equation*}
Where $P$ is probability measure and $p^*$ - density of this measure (if it exists).
For equation~\ref{rel sde} we can obtain:
\begin{equation*}
    \mathcal{A}_{\alpha} \phi = L_{-\nabla{f}} \phi + \frac{\alpha^2}{2} \|\nabla{f}\|_2^2 \Laplace{\phi}.
\end{equation*}
We will look for a stationary solution to the Fokker-Planck equation for the equation~\ref{rel sde}:
\begin{equation*}
    \forall \phi \in C_0^{\infty}\left(\mathbb{R}^n\right) \int_{\mathbb{R}^n} \left(\mathcal{A}_{\alpha} \phi \right) P^*(dx) = 0.
\end{equation*}
We can note, that any distribution $P^*$ concentrated on $X^* = \lbrace x \; | \; \nabla{f}(x) = 0 \rbrace$ will be solution. Thus, assuming a model~\ref{relative_inexact}, we can expect same convergence as for models without noise. Note, that~\ref{rel sde} does not assume $\alpha < 1$, this effect will be noted in Section~\ref{section stm}.

\section{Gradient descent} \label{section stoh gd correlated}
In this section, we consider problem~\ref{optim}. We will use conditions~\ref{abs_inexact} and~\ref{correl}. Let us introduce classic procedure
\begin{algorithm}[H]
\caption{Gradient Descent $(h, x^0)$}
\label{alg gd}
\begin{algorithmic}[1]
\State 
\noindent {\bf Input:} Starting point $x^0$, number of steps $N$, $h$ - learning rate.
\For {$k = 1 \dots N$}
        \State $x^{k} = x^{k - 1} - h \widetilde{\nabla}f(x^{k - 1})$
\EndFor
\State 
\noindent {\bf Output:} $x^N$.
\end{algorithmic}
\end{algorithm}
\begin{lemma} \label{gd step}
    Let $f$ function satisfies condition~\ref{smooth cond}, and $\widetilde{\nabla}f$ satisfies~\ref{growth condition} and~\ref{correl}. Then Algorithm~\ref{alg gd} with $h < \frac{2\nu\gamma}{L \rho^2}$ produce:
    \begin{equation*}
        \|\nabla f(x^k) \|_2^2 \leqslant \tau \left( f(x^{k}) - f(x^{k + 1}) \right)
    \end{equation*}
    Where:
    \begin{equation*}
        \tau^{-1} = h \left(\nu \gamma - \frac{Lh\rho^2}{2} \right)
    \end{equation*}
\end{lemma}
\begin{proof}
    \begin{equation*}
        \begin{gathered}
            f(x^{k + 1}) \leqslant f(x^k) + \langle \nabla{f}(x^k), x^{k + 1} - x^k \rangle + \frac{L}{2} \|x^{k + 1} - x^k \|_2^2, \\
            f(x^{k + 1}) \leqslant f(x^k) + \langle \nabla{f}(x^k), -h\widetilde{\nabla}f(x^k) \rangle + \frac{h^2L}{2} \| \widetilde{\nabla} f(x^k) \|_2^2, \\
            f(x^{k + 1}) \leqslant f(x^k) - h \nu \gamma \| \nabla{f}(x^k) \|_2^2 + \frac{h^2L \rho^2}{2} \| \nabla{f}(x^k) \|_2^2
        \end{gathered} 
    \end{equation*}
Then we obtain desired inequality.
\end{proof}
\begin{lemma} \label{gd dist bound}
    Let $f$ function satisfies conditions Lemma~\ref{gd step}, then:
    \begin{equation*}
        \|x^k - x^* \|_2^2 \leqslant 2 \left(L h^2 \rho^2 \tau + 1\right) R^2
    \end{equation*}
\end{lemma}
\begin{proof}
    Using Lemma~\ref{gd step} and telescopic sum:
    \begin{equation*}
        \begin{gathered}
            \|x^k - x^* \|_2^2 \leqslant 2\|x^k - x^0 \|_2^2 + 2\|x^0 - x^* \|_2^2, \\
            \|x^k - x^0 \|_2^2 \leqslant \| h \sum_{j = 0}^{k - 1} \widetilde{\nabla}f(x^j) \|_2^2 \leqslant 2 h^2 \rho^2 \sum_{j = 0}^{k - 1} \| \nabla f(x^j) \|_2^2 \leqslant
            2 h^2 \rho^2 \tau \Delta_0 \leqslant L h^2 \rho^2 \tau R^2, \\
            \|x^k - x^* \|_2^2 \leqslant 2 \left(L h^2 \rho^2 \tau + 1\right) R^2
        \end{gathered} 
    \end{equation*}
\end{proof}
\begin{theorem} \label{gd convex}
    Let $f$ - convex function and satisfies condition~\ref{smooth cond}, and $\widetilde{\nabla}f$ satisfies~\ref{growth condition} and~\ref{correl} then Algorithm~\ref{alg gd}. Then with $h < \frac{2\nu\gamma}{L \rho^2}$ we obtain:
    \begin{equation*}
        \begin{gathered}   
        f(x^N) - f^* \leqslant \frac{2 \left(L h^2 \rho^2 \tau + 1\right)}{h \left(\nu \gamma - \frac{Lh\rho^2}{2} \right)} \frac{R^2}{N + 1}, \\
        R = \|x^0 - x^* \|_2, \; \tau^{-1} = h \left(\nu \gamma - \frac{Lh\rho^2}{2} \right)
        \end{gathered}
    \end{equation*}
\end{theorem}
\begin{proof}
    Let us define: $\Delta_k = f(x^k) - f^*$.
    Using convexity of function $f$ and Cauchy Bunyakovsky Schwarz inequality we obtain:
    \begin{equation*}
        \begin{gathered}
            f(x^k) + \langle \nabla{f}(x^k), x^* - x^k \rangle \leqslant f^*, \\
            \Delta_{k} \leqslant \|\nabla f(x^k) \|_2 \|x^k - x^* \|_2 
        \end{gathered}
    \end{equation*}
    Let us define $\lambda = 2 \left(L h^2 \rho^2 \tau + 1\right)$. Using Lemma~\ref{gd step} and Lemma~\ref{gd dist bound}:
    \begin{equation*}
        \begin{gathered}
            \| \nabla f(x^k) \|_2^2 \geqslant \frac{\Delta_k}{\lambda R^2}, \\
            \Delta_{k + 1} \leqslant \Delta_k - \frac{1}{\lambda \tau R^2} \Delta_k^2. 
        \end{gathered}
    \end{equation*}
    Then, using previous relation we will prove by induction convergence rate:
    Base $N = 2$, because $\tau$ minimum takes place with $h = \frac{\nu \gamma}{L \rho^2}$:
    \begin{equation*}
       \begin{gathered}
           f(x^2) \leqslant f^* + \frac{L}{2} \|x^2 - x^* \|_2^2, \\
           \Delta_2 \leqslant \frac{\lambda L R^2}{2} < \frac{2 \lambda L R^2}{3} \leqslant \frac{\lambda \tau R^2}{3}
       \end{gathered} 
    \end{equation*}
    Step $N + 1 > 2$:
    \begin{equation*}
       \begin{gathered}
            \text{Assume: } \Delta_N \leqslant \frac{\lambda \tau R^2}{N + 1}, \\
            \text{If: } \Delta_{N} \leqslant \frac{\lambda \tau R^2}{N + 2} \Rightarrow \Delta_{N + 1} \leqslant \frac{\lambda \tau R^2}{N + 2}, \\
            \text{If: } \Delta_N > \frac{\lambda \tau R^2}{N + 2}, \\
            \text{Function: } g(x) = \frac{1}{x} - \frac{1}{x^2} \text{ - decreasing for } x > 2, \\
            \frac{\lambda \tau R^2}{N + 1} - \frac{\lambda \tau R^2}{(N + 1)^2} < \frac{\lambda \tau R^2}{N + 2} \Rightarrow \Delta_{N + 1} \leqslant \frac{\lambda \tau R^2}{N + 2}
       \end{gathered} 
    \end{equation*}
    Thus we obtain convergence.
\end{proof}
\begin{theorem} \label{gd nonconvex}
    Let $f$ - function and satisfies condition~\ref{smooth cond}, $\widetilde{\nabla}f$ satisfies~\ref{growth condition} and~\ref{correl} then procedure~\ref{alg gd}. Then with $h < \frac{2\nu\gamma}{L \rho^2}$ we obtain:
    \begin{equation*}
        \min_{0 \leqslant k \leqslant N} \| \nabla{f}(x^k) \|_2^2 \leqslant
        \frac{f(x^0) - f^*}{h \left(\nu \gamma - \frac{Lh\rho^2}{2} \right) \left(N + 1\right)}.
    \end{equation*}
\end{theorem}
\begin{proof}
    Using Lemma~\ref{gd step} and summarize it from $0$ to $N$:
    \begin{equation*}
        \begin{gathered}
            \sum_{j = 0}^N \| \nabla{f}(x^j) \|_2^2 \leqslant \tau \sum_{j = 0}^N \left( f(x^{j}) - f(x^{j + 1}) \right) = \tau\left(f(x^0) -  f(x^{N + 1}) \right) \leqslant \tau \left( f(x^0) - f^* \right), \\
            (N + 1) \min_{0 \leqslant k \leqslant N} \| \nabla{f}(x^k) \|_2^2 \leqslant \tau \left( f(x^0) - f^* \right), \\
            \min_{0 \leqslant k \leqslant N} \| \nabla{f}(x^k) \|_2^2 \leqslant \frac{\tau \left(f(x^0) - f^* \right)}{N + 1}.
        \end{gathered}
    \end{equation*}
\end{proof}

\begin{remark} \label{h choice gd}
    In conditions of Theorem~\ref{gd convex} we can choice
    \begin{equation*}
        h = \frac{\nu \gamma}{L \rho^2}.
    \end{equation*}
    Then:
    \begin{equation*}
        \begin{gathered}
            \tau = \frac{2L \rho^2}{\nu^2 \gamma^2}, \\
            \lambda = 2 \left(L h^2 \rho^2 \tau + 1\right) = 6, \\
            f(x^N) - f^* \leqslant \frac{12 L \rho^2 R^2}{\nu^2 \gamma^2(N + 1)} = O \left( \frac{\rho^2}{\nu^2 \gamma^2} \frac{LR^2}{N + 1} \right), \\
            \min_{0 \leqslant k \leqslant N} \| \nabla{f}(x^k) \|_2^2 \leqslant
            \frac{\rho^2}{\nu^2 \gamma^2} \frac{2L \left(f(x^0) - f^* \right)}{N + 1}.
        \end{gathered}
    \end{equation*}
\end{remark}
\begin{theorem} \label{alpha for gd}
    Let $f$ - function and satisfies condition~\ref{smooth cond}, and $\widetilde{\nabla}f$ satisfies~\ref{relative_inexact} then Algorithm~\ref{alg gd} with 
    \begin{equation*}
        h = \left(\frac{1 - \alpha}{1 + \alpha} \right)^{\frac{3}{2}} \frac{1}{L},
    \end{equation*}
    has convergence:
    \begin{equation*}
        \begin{gathered}
            f(x^N) - f^* \leqslant \frac{12 L (1 + \alpha) R^2}{(1 - \alpha)^3 (N + 1)} = O \left( \frac{1 + \alpha}{(1 - \alpha)^3} \frac{LR^2}{N + 1} \right), \\
            \min_{0 \leqslant k \leqslant N} \| \nabla{f}(x^k) \|_2^2 \leqslant
            \frac{1 + \alpha}{(1 - \alpha)^3}\frac{2L \left(f(x^0) - f^* \right)}{N + 1}.
        \end{gathered}
    \end{equation*}
\end{theorem}
\begin{proof}
    Continuing Remark~\ref{h choice gd} and note~\ref{alpha to growth} we obtain theorem above.
\end{proof}
\begin{remark} \label{lyapunov gd proof}
    We provide not default proof for gradient descent convergence. The reason for this may be the Lyapunov function tool for a dynamic system:
    \begin{equation*}
        \dot{x} = -\widetilde{\nabla}f(x).
    \end{equation*}
    Function $f$ - will be Lyapunov function for system above, if $\widetilde{\nabla} f$ satisfies condition~\ref{correl}. However $\| x - x^* \|_2^2$ will not:
    \begin{equation*}
        \begin{gathered}
            L_{\widetilde{\nabla}f} f = -\langle \widetilde{\nabla} f(x), \nabla{f}(x) \rangle \leqslant - \gamma \|\widetilde{\nabla} f(x) \|_2 \|\nabla{f}(x) \|_2 < 0, \\
            L_{\widetilde{\nabla}f} \| x - x^* \|_2^2 = - \langle \widetilde{\nabla} f(x), x - x^* \rangle \Rightarrow \text{sign of this product is ambiguous }.
        \end{gathered}
    \end{equation*}
\end{remark}

\section{Similar triangle methods and stochastic optimization} \label{section stm}

For constrained optimization with set $Q$ we can provide such method:

\begin{algorithm}[H]
\caption{STM $(L,  x_{start})$, $Q \subseteq \mathbb{R}^n$}
	\label{stm}
\begin{algorithmic}[1]
\State 
\noindent {\bf Input:} Starting point $x_{start}$, number of steps $N$.
\State {\bf Set} $\tilde x^0 = x_{start}$, $\alpha_0 = \frac{1}{L}$, $A_0 = \frac{1}{L}$.
\State {\bf Set} $\psi_0(x) = \frac{1}{2} \|x - \tilde{x}^0\|_2^2 + \alpha_0\left(f(\tilde{x}^0) + \langle \tilde{\nabla} f(\tilde{x}^0), x - \tilde{x}^0 \rangle  \right)$. 
\State {\bf Set} $z_0 = \argmin_{y \in Q}{\psi_0(y)}$, $x^0 = z^0$.
\For {$k = 1 \dots N$}
        \State Find $\alpha_k$ from $(1 + A_{k - 1})(A_{k - 1} + \alpha_k) = L \alpha_k^2$, \label{step:alpha_k}\\
        \State or equivalently $\alpha_k = \frac{1}{2L} + \sqrt{\frac{1}{4L^2} + \frac{A_{k - 1}}{L}} $, 
        \State $A_k = A_{k - 1} + \alpha_k,$
        \State $\tilde x^k = \frac{A_{k - 1} x^{k - 1} + \alpha_k z^{k - 1}}{A_k},$
        \State $\psi_k(x) = \psi_{k - 1}(x) + \alpha_k \left(f(\tilde{x}^k) + \langle \widetilde{\nabla} f(\tilde{x}^k), x - \tilde{x}^k \rangle \right),$
        \State $z^k = \argmin_{y \in Q} \psi_k(y),$ 
        \State $x^k = \frac{A_{k - 1} x^{k - 1} + \alpha_k z^k}{A_k}.$
    %\end{eqnarray*}
    %\end{align*}
\EndFor
\State 
\noindent {\bf Output:} $x_N$.
\end{algorithmic}
\end{algorithm}

For not constrained optimization we can rewrite it as implicit method.

\begin{algorithm}[H]
\caption{STM (implicit) $(L,  x_{start})$}
	\label{alg stm rn}
\begin{algorithmic}[1]
\State 
\noindent {\bf Input:} Starting point $x_{start}$, number of steps $N$.
\State {\bf Set} $y^0 = x_{start}$, $\alpha_0 = \frac{1}{L}$, $A_0 = \frac{1}{L}$.
\State {\bf Set} $z_0 = y^0 - \alpha_0 \widetilde{\nabla} f(y^0)$, $x^0 = z^0$.
\For {$k = 1 \dots N$}
        \State $\alpha_k = \frac{1}{2L} + \sqrt{\frac{1}{4L^2} + \frac{A_{k - 1}}{L}} $, 
        \State $A_k = A_{k - 1} + \alpha_k,$
        \State $y^k = \frac{A_{k - 1} x^{k - 1} + \alpha_k z^{k - 1}}{A_k},$
        \State $z^k = z^{k - 1} - \alpha_k \widetilde{\nabla} f(y^{k})$ 
        \State $x^k = y^{k} - \frac{1}{L} \widetilde{\nabla} f(y^{k})$
\EndFor
\State 
\noindent {\bf Output:} $x_N$.
\end{algorithmic}
\end{algorithm}
\subsection{Stochastic optimization}
In this section we will consider stochastic optimization:
\begin{equation*}
    \mathbb{E}_{\xi} \widetilde{\nabla} f(x, \xi) = \nabla{f}(x).
\end{equation*}
With analogue growth condition:
\begin{equation*}
    \mathbb{E}_{\xi} \|\widetilde{\nabla} f(x, \xi) \|_2^2 \leqslant \kappa \| \nabla f(x) \|_2^2.
\end{equation*}
In paper~\cite{vaswani2019fast}, Theorem 2 was proved about the convergence of a method similar to ~\ref{alg stm rn} with step correction:
\begin{equation*}
    \mathbb{E} f(x^N) - f^* \leqslant \frac{2\kappa^2 L}{N^2} \|x^N - x^* \|_2^2.
\end{equation*}
The condition of unbiased gradient turned out to be very important in the proof of this theorem. First of all it is easy to get lower bound of growth condition~\ref{growth condition}:
\begin{equation*}
    \mathbb{E} \|\widetilde{\nabla} f(x, \xi) \|_2 \geqslant \|\mathbb{E} \widetilde{\nabla} f(x, \xi) \|_2 = \| \nabla{f}(x) \|_2.
\end{equation*}
Gradient direction condition~\ref{correl} can be interpreted a little differently. We assume that $v$ corresponds to the coordinate $x_1$.
Define $\mathbb{S}^{n - 1}_{v+}$ - top hemisphere for direction $v$. Let us define simple noised gradient model for analysis:
\begin{equation*}
    \begin{gathered}
        \mathbb{E} u = v, \\
        u \sim \omega \mathcal{U}\left(\mathbb{S}^{n - 1}_{v+} \right) + (1 - \omega) \mathcal{U}\left(\mathbb{S}^{n - 1}_{v-} \right), \; \omega \in [0; 1].
    \end{gathered}
\end{equation*}
Then we can calculate norm of vector $v$, it can be done, because it has only $x_1$ non zero coordinate:
\begin{equation*}
    \begin{gathered}
        g(n) = \frac{1}{\text{area}\left(\mathbb{S}^{n - 1}_{v+} \right)} \int_{\mathbb{S}^{n - 1}_{v+}} x_1 dA, \\
        \|v \|_2 = (2\omega - 1) g(n) \Rightarrow \omega = \frac{1}{2} + \frac{\|v\|_2}{2 g(n)} > \frac{1}{2}, \\
        P\left(\langle u, v \rangle > 0 \right) = \omega > \frac{1}{2}.
    \end{gathered}
\end{equation*}
Then we can move to batching estimation to study direction of batched estimation:
\begin{equation*} \label{model uni}
    \begin{gathered}    
        S_m = \frac{1}{m} u_j, \; u_j - \text{i.i.d (batching)}, \\
        P \left( \langle S_m, v \rangle > 0 \right) = P \left( \frac{1}{m}\sum_{k = 1}^m \cos{\phi_k} > 0 \right), \; \phi_k \sim \mathcal{D}_{\omega}, \\
         \mathcal{D}_{\omega} = \omega \mathcal{U}\left(\left[0; \frac{\pi}{2}\right] \right) + (1 - \omega) \mathcal{U}\left(\left[\frac{\pi}{2}; \pi \right] \right), \\
         a_{\omega} = \mathbb{E}_{\mathcal{D}_{\omega}} \cos{\phi} = \frac{2(2\omega - 1)}{\pi}, \\
         \sigma_{\omega}^2 = \Var_{\mathcal{D}_{\omega}} \sin{\phi} = \frac{1}{2} - a_{\omega}^2, \\
         P\left( \langle S_m, v \rangle > 0 \right) \geqslant 1 - P\left( \left| \frac{1}{m}\sum_{k = 1}^m \cos{\phi_k} - a \right| \geqslant a \right) \geqslant 1 - \frac{1}{m} \left(\frac{1}{2a_{\omega}^2} - 1\right), \\
         P\left( \langle S_m, v \rangle > 0 \right) \geqslant 1 - \frac{1}{m} \left(\frac{\pi^2}{8(2\omega - 1)^2} - 1 \right) = 1 - \frac{1}{m} \left(\frac{\pi^2}{8} \frac{g^2(n)}{\|v\|_2^2} - 1 \right).
    \end{gathered}
\end{equation*}
Such model gives us, that unbiased uniformly distributed $\widetilde{\nabla}f(x, \xi)$ gives $P\left(\langle \widetilde{\nabla}f(x, \xi), \nabla f(x) \rangle > 0 \right) = \omega > \frac{1}{2}$ and using batching this probability can be increased to $1 - \frac{1}{m}\left(\frac{\pi^2}{8} \frac{1}{(2 \omega - 1)^2} - 1 \right)$ and we can note as the $\nabla{f}(x)$ norm increases, so does the probability.
\subsection{Strongly convex case, not stochastic} 
We call function $f$ $\mu$-strongly convex, if:
\begin{equation} \label{strong mu}
    (\forall x, y \in \mathbb{R}^n): \; f(x) + \langle \nabla{f}(x), y - x \rangle + \frac{\mu}{2} \|x -  y \|_2^2 \leqslant f(y)
\end{equation}
In \cite{kornilov2023intermediate} convergence of analogue of Algorithm~\ref{alg stm rn} for strongly convex case~\ref{strong mu} using restarts:
\begin{equation*}
    \begin{gathered}
        N = O \left( \sqrt{\frac{L}{\mu}} \log{\frac{\mu R_0^2}{\varepsilon}} \right) \Rightarrow f(x^N) - f^* \leqslant \varepsilon, \\
        \text{If } \alpha \text{(relative error level)} = O \left( \sqrt{\frac{\mu}{L}} \right)
    \end{gathered}
\end{equation*}
We can give motivation for $\sqrt{\frac{\mu}{L}}$ bound, which appeared in \cite{kornilov2023intermediate, gannot2021frequency} proofs. Firstly one solid geometry fact:
\begin{lemma}[Cosine $\mathbb{R}^3$ theorem] \label{cosinus theorem}
    Let $OABC$ - tetrahedron, $\alpha = \angle BOC, \beta = \angle AOC, \gamma = \angle AOB, \delta = \angle(AOB, AOC)$, then
    \begin{equation*}
        \cos \alpha = \cos \beta \cos \gamma + \sin \beta \sin \gamma \cos \delta.
    \end{equation*}
\end{lemma}
\begin{lemma}
    We call function $f$ - noised $\beta, \mu$-strongly quasiconvex at points $x, x^*$, if:
    \begin{equation} \label{squasi convex}
        f(x) + \beta^{-1} \langle \widetilde{\nabla}f(x), x^* - x \rangle + \frac{\mu}{2} \|x^* - x \|_2^2 \leqslant f^*, \; \beta \in (0; 1). 
    \end{equation}
    If $f$ - $\mu$-strongly convex~\ref{strong mu} and smooth~\ref{smooth cond}, $\widetilde{\nabla}f$ satisfies~\ref{relative_inexact}, $\alpha < \sqrt{\frac{\mu}{L}}$, then $f$ satisfies~\ref{squasi convex} with $\beta = \frac{1}{2}$.
\end{lemma}
\begin{proof}
    As $f$ is convex, then:
    \begin{equation*}
        f(x) + \langle \nabla{f}(x), x^* - x \rangle \leqslant f^* \Rightarrow \langle \nabla{f}(x), x - x^* \rangle \geqslant f(x) - f^*.
    \end{equation*}
    Using $\mu, L$ constants conditions:
    \begin{equation*}
        \langle \nabla{f}(x), x - x^* \rangle = \gamma(x) \|\nabla{f}(x) \|_2 \|x^* - x \|_2 \leqslant 2\gamma(x) \sqrt{\frac{L}{\mu}} \left(f(x) - f^* \right).
    \end{equation*}
    Then we get $\gamma(x) \geqslant \frac{1}{2} \sqrt{\frac{\mu}{L}}$. We can decompose:
    \begin{equation*}
        \widetilde{\nabla} f(x) = \nabla f(x) + r(x), \; \|r(x) \|_2 \leqslant \alpha \| \nabla f(x) \|_2.
    \end{equation*}
    Let us define $\widetilde{\gamma}(x)$ such as:
    \begin{equation*}
        \left\langle \frac{1}{2}\nabla f(x) + r(x), x - x^* \right\rangle = \widetilde{\gamma}(x) \left\|\frac{1}{2}\nabla f(x) + r(x) \right\|_2 \| x - x^* \|_2.
    \end{equation*}
    Using Lemma~\ref{cosinus theorem}:
    \begin{equation*}
        \widetilde{\gamma}(x) \geqslant \frac{1}{2} \sqrt{1 - \frac{\alpha^2}{4}} \sqrt{\frac{\mu}{L}} - \frac{\alpha}{2} \sqrt{1 - \frac{\mu}{4L}}.
    \end{equation*}
    That is, $\left\langle \frac{1}{2} \nabla{f}(x) + r(x), x - x^* \right\rangle \geqslant 0 \Longleftrightarrow \alpha < \sqrt{\frac{\mu}{L}}$. Then:
    \begin{align*}
        &  f(x) + 2 \left\langle \nabla{f}(x) + r(x), x^* - x \right\rangle + \frac{\mu}{2}\|x^* - x \|_2^2 \\
        & = f(x) + \langle \nabla{f}(x), x^* - x \rangle + \frac{\mu}{2} \| x^* - x \|_2^2 + 2\left\langle \frac{1}{2}\nabla{f}(x) + r(x), x^* - x \right\rangle \\
        & \leqslant f^*.
    \end{align*}
\end{proof}

\section{Conjugate Gradients} \label{CG}

In this section we will provide accelerated method~$\ref{alg cg}$ for model~\ref{smooth cond}, \ref{correl}, \ref{growth condition}. We will limit the set of functions to the simplified strong convex condition. We will say, that function $f$ satisfies quadratic growth condition with parameter $\mu$ at point $x$ if:
\begin{equation} \label{mu cond}
    \frac{\mu}{2} \|x - x^* \|_2^2 \leqslant f(x) - f^*.
\end{equation}
This condition was taken from~\cite{guminov2023accelerated} as well as the algorithm~\ref{alg cg}. We will use this condition to prove analogues of the results from the papers~\cite{gannot2021frequency, kornilov2023intermediate, vasin2023accelerated}.

\begin{algorithm}[H]
\caption{CG $(h,  x_{start})$}
	\label{alg cg}
\begin{algorithmic}[1]
\State 
\noindent {\bf Input:} Starting point $x_{start}$, number of steps $N$.
\State {\bf Set} $u^0 = 0$.
\State {\bf Set} $y^0 = 0$.
\For {$k = 1 \dots N$}
        \State $x^k = y^{k - 1} - h \widetilde{\nabla}f(y^{k - 1})$
        \State $u^k = u^{k - 1} + \widetilde{\nabla}f(y^{k - 1})$
        \State $y^k = \argmin\limits_{x = x^0 + \text{Lin}\lbrace{ x^k - x^0, u^k \rbrace}}f(x).$
\EndFor
\State 
\noindent {\bf Output:} $x^N$.
\end{algorithmic}
\end{algorithm}

\begin{theorem} \label{cg conv}
    Let $f$ - convex function and satisfies condition~\ref{smooth cond}, $\widetilde{\nabla}f$ satisfies~\ref{correl} and~\ref{growth condition}, starting point $x^0$ satisfies~\ref{mu cond}, then procedure~\ref{alg gd} with $h < \frac{2\nu\gamma}{L \rho^2}$ produce:
    \begin{equation*}
        f(x^N) - f^* \leqslant \omega (f(x^0) - f^*), \; \omega = \frac{3}{4}.
    \end{equation*}
    Where:
    \begin{equation*}
        N > \max \left\lbrace \frac{2}{\omega} \rho \sqrt{\frac{\tau}{\mu}}, \frac{4 \rho^2 \sqrt{1 - \gamma^2} \tau}{\mu \omega^2} \right\rbrace, \; \tau^{-1} = h \left(\nu \gamma - \frac{Lh\rho^2}{2} \right)
    \end{equation*}
\end{theorem}
\begin{proof}
    Let's conduct a proof against the contrary, assume:
    \begin{equation*}
        f(x^N) - f^* \geqslant \omega (f(x^0) - f^*).
    \end{equation*}
    Using Lemma~\ref{gd step} and $f(x^k) \leqslant f(y^k)$ we obtain:
    \begin{equation*}
        \| \nabla{f}(y^k) \|_2^2 \leqslant \tau \left( f(x^k) - f(x^{k + 1}) \right)
    \end{equation*}
Summing up inequalities above:
\begin{equation*}
    \sum_{k = 0}^{N - 1} \|\nabla{f}(y^k) \|_2^2 \leqslant \tau \left(f(x^0) - f^* \right).    
\end{equation*}
From definition of $y^k$ and first order condition we obtain:
\begin{equation*}
    \langle \nabla{f}(y^k), x^k - x^0 \rangle = \langle \nabla{f}(y^k), u^k \rangle = 0.
\end{equation*}
Using convexity:
\begin{equation*}
    \omega (f(x^0) - f^*) \leqslant f(y^k) - f^* \leqslant \langle \nabla{f}(y^k), x^* - y^k \rangle = \langle \nabla{f}(y^k), x^* - x^0 \rangle.
\end{equation*}
Summing up inequalities above:
\begin{equation*}
    \omega N (f(x^0) - f^*) \leqslant \langle u^N, x^* - x^0 \rangle.
\end{equation*}
Using condition~\ref{correl}, Lemma~\ref{cosinus theorem} and $\langle u^k, \nabla{f}(y^k) \rangle = 0$:
\begin{equation*}
    \langle \widetilde\nabla{f}(y^k), u^k \rangle \leqslant \sqrt{1 - \gamma^2} \|\widetilde\nabla{f}(y^k) \|_2 \|u^k \|_2.
\end{equation*}
Then:
\begin{align*}
    \left \| u^N \right \|_2^2
    & \leqslant \| \widetilde{\nabla}{f}(y^{N-1}) \|_2^2 + \|u^{N - 1} \|_2^2 + \sqrt{1 - \gamma^2} \| u^{N-1}\|_2 \|\widetilde{\nabla}{f}(y^{N-1})\|_2 \\
    & \leqslant \rho^2 \| \nabla{f}(y^{N-1}) \|_2^2 + \|u^{N - 1} \|_2^2 + \sum_{k = 0}^{N - 2} \rho^2 \sqrt{1 - \gamma^2} \left(\| \nabla{f}(y^k) \|_2^2 + \|\nabla{f}(y^{N-1}\|_2^2\right) \\
    & \leqslant \rho^2 \sum_{k = 0}^{N - 1} \| \nabla{f}(y^{k}) \|_2^2 + \rho^2 \sqrt{1 - \gamma^2} \sum_{k = 0}^{N - 1} (N - 1) \| \nabla{f}(y^k) \|_2^2 \\ 
    & \leqslant \rho^2 \left(1 + \sqrt{1 - \gamma^2} N \right) \sum_{k = 0}^{N - 1} \|\nabla{f}(y^k) \|_2^2 \leqslant 2\rho^2 \max\lbrace1, \sqrt{1 - \gamma^2} N\rbrace \tau (f(x^0) - f^*)
\end{align*}
Using quadratic growth at point $x^0$:
\begin{equation*}
    \frac{\mu}{2} \|x^0 - x^* \|_2^2 \leqslant f(x^0) - f(x^*).
\end{equation*}
Then:
\begin{equation*}
    \begin{gathered}
        \omega N (f(x^0) - f^*) \leqslant \sqrt{\frac{4}{\mu}} \rho \sqrt{\tau \max\lbrace 1, \sqrt{1 - \gamma^2}N \rbrace } (f(x^0) - f^*) \Rightarrow \\
        \Rightarrow N \leqslant \max \left\lbrace \frac{2}{\omega} \rho \sqrt{\frac{\tau}{\mu}}, \frac{4 \rho^2 \sqrt{1 - \gamma^2} \tau}{\mu \omega^2} \right\rbrace
    \end{gathered}
\end{equation*}
This ends the proof.
\end{proof}

\begin{theorem} \label{cg alpha}
    Let $f$ - function and satisfies condition~\ref{smooth cond}, $\widetilde{\nabla}f$ satisfies~\ref{relative_inexact} and $x^0$ satisfies~\ref{mu cond} then procedure~\ref{alg cg} with 
    \begin{equation*}
        h = \left(\frac{1 - \alpha}{1 + \alpha} \right)^{\frac{3}{2}} \frac{1}{L},
    \end{equation*}
    and:
    \begin{equation*}
        \alpha < \sqrt{\frac{\mu}{L}}.
    \end{equation*}
    produce:
    \begin{equation*}
        \begin{gathered}
            f(x^N) -f^* \leqslant \frac{3}{4} (f(x^0) - f^*), \\
            N > \frac{64 \left(1 + \alpha\right)^{3}}{9\left(1 - \alpha \right)^3} \sqrt{\frac{L}{\mu}}.
        \end{gathered}
    \end{equation*}
\end{theorem}
\begin{proof}
    Using Theorem~\ref{cg conv} and defined step $h$ we obtain the result.
    We should note, that
    \begin{equation*}
        \frac{\left(1 + \alpha \right)^{3}}{\left(1 - \alpha \right)^3} = O\left( 1 \right).
    \end{equation*}
\end{proof}
Theorem~\ref{cg alpha} allows us to obtain an accelerated method that preserves convergence under conditions of relative noise with $\alpha = O\left(\sqrt{\frac{\mu}{L}} \right)$. This result corresponds papers~\cite{gannot2021frequency, kornilov2023intermediate}. The disadvantage of this approach is the need to use the oracle of low-dimensional optimization.

\section{Sufficient conditions for datasets} \label{sufficient dataset}
\subsection{Logistic regression}
In this section we will provide some sufficient conditions for~\ref{correl},~\ref{growth condition}. Consider sum-structured optimization problem~\ref{sum opt} with dataset $\mathcal{D} = \lbrace (x_k, y_k) \rbrace_{k = 1}^M$, where $x_k$ - features of object $k$ and $y_k \in \lbrace 0, 1 \rbrace$ - object`s class, for model we consider logistic regression for simplicity:
\begin{equation} \label{logistic}
    \begin{gathered}
        \sigma(x) = \frac{1}{1 + e^{-x}}, \\
        f_k(\theta) = y_k \ln \sigma\left(x_k^T \theta \right) + (1 - y_k) \ln \left( 1 - \sigma\left(x_k^T \theta \right) \right), \\
        f(\theta) = \frac{1}{M} \sum_{j = 1}^M f_k(\theta), \\
        u_k(\theta) = \nabla f_k(\theta) = x_k \left( y_k \left( 1 - \sigma\left(x_k^T \theta \right) \right) - (1 - y_k) \sigma\left(x_k^T \theta \right) \right).
    \end{gathered}
\end{equation}
We define features consistency property for dataset $\mathcal{D}$:
\begin{equation} \label{consist data}
    \begin{gathered}
        (\exists \Upsilon \in (0; 1]) (\forall i, k \in \lbrace 1 \dots M \rbrace), \\ 
        y_i = y_k \Rightarrow \langle x_i, x_k \rangle \geqslant \Upsilon \| x_i \|_2 \|x_k\|_2, \\
        y_i \not= y_k \Rightarrow \langle x_i, x_k \rangle \leqslant -\Upsilon \| x_i \|_2 \|x_k\|_2, \\
    \end{gathered}
\end{equation}
\begin{theorem} \label{gamma cons}
    If dataset $\mathcal{D}$ satisfies features consistency~\ref{consist data} then single element gradient estimation satisfies condition~\ref{correl}.
\end{theorem}
\begin{proof}
    Let $\widetilde{\nabla}f(x) = u_{\xi}(x)$, $\xi \in \lbrace 1 \dots M \rbrace$.
    \begin{equation*}
        \begin{gathered}
            \langle u_{\xi}(\theta), \frac{1}{M} \sum_{j = 1}^M u_k(\theta) \rangle = \frac{1}{M} \sum_{k = 1}^M \langle u_\xi(\theta), u_k(\theta) \rangle, \\
            y_{\xi} = y_k = 1 \Rightarrow \langle u_\xi(\theta), u_k(\theta) \rangle = \langle x_{\xi}, x_k \rangle \left( 1 - \sigma\left(x_{\xi}^T \theta \right) \right) \left( 1 - \sigma\left(x_k^T \theta \right) \right) , \\
            y_{\xi} = 1, y_k = 0 \Rightarrow \langle u_\xi(\theta), u_k(\theta) \rangle = -\langle x_{\xi}, x_k \rangle \left( 1 - \sigma\left(x_{\xi}^T \theta \right) \right) \sigma\left(x_k^T \theta \right) \geqslant \Upsilon \| u_{\xi}(\theta) \|_2 \| u_k(\theta) \|_2, \\
            \text{Analogically for } y_{\xi} = y_k = 0, \\
            \frac{1}{M} \sum_{k = 1}^M \langle u_\xi(\theta), u_k(\theta) \rangle \geqslant \Upsilon \| u_{\xi}(\theta) \|_2  \bigg \| \frac{1}{M} \sum_{k = 1}^M u_k(\theta) \bigg \|_2, \\
            \Rightarrow \langle \widetilde{\nabla}f(x), \nabla f(x) \rangle \geqslant \Upsilon \| \widetilde{\nabla}f(x) \|_2 \| \nabla f(x) \|_2.
        \end{gathered}
    \end{equation*}
\end{proof}
\begin{theorem}
    If dataset $\mathcal{D}$ satisfies features consistency~\ref{consist data} then batch gradient estimation satisfies condition~\ref{correl}.
\end{theorem}
\begin{proof}
    Using linearity of scalar product and Theorem~\ref{gamma cons} easy to obtain theorem.
\end{proof}
\subsection{Linear inverse problem}
Consider a system of linear equations with respect to $\theta$:
\begin{equation} \label{linear eq}
    X\theta = Y, \; X \in \mathrm{GL}_m(\mathbb{R}) - \text{invertable matrix}, \; Y \in \mathbb{R}^m.
\end{equation}
We can map the convex optimization problem to a system~\ref{linear eq}:
\begin{equation} \label{lin inv opt}
    f(\theta) = \frac{1}{2m}\|X \theta - Y \|_2^2 = \frac{1}{2m} \sum_{k = 1}^m (x_k^T \theta - y_k)^2.
\end{equation}
Since the matrix is non-generated, then $\exists \theta^*: \; X\theta^* = Y$, then:
\begin{equation*}
    \begin{gathered}
        \nabla f(\theta) = \frac{1}{m} \sum_{k = 1}^m P_k \Delta \theta, \\
        \Delta \theta = \theta - \theta^*, \\
        P_k = x_k x_k^T, \\
        S = \frac{1}{m} \sum_{k = 1}^m P_k.
    \end{gathered}
\end{equation*}
\begin{theorem}
    If matrix $X \in \mathrm{GL}_m(\mathbb{R})$ and satisfies:
    \begin{equation*} \label{rows cond}
        \forall j, k: \; 1 \leqslant j, k \leqslant m: \; \langle x_k, x_j \rangle \geqslant \Upsilon \|x_k\|_2 \|x_j\|_2.
    \end{equation*}
    Then single element gradient estimation $\widetilde{\nabla} f(\theta) = P_{\xi} \Delta \theta$ satisfies~\ref{correl}.
\end{theorem}
\begin{proof}
    \begin{align*}
        \langle \nabla f(\theta), \widetilde{\nabla} f(\theta) \rangle & = \langle S \Delta \theta, P_{\xi} \Delta \theta \rangle = \frac{1}{m} \sum_{k = 1}^m \langle P_k \Delta \theta, P_{\xi} \Delta \theta \rangle \geqslant \\
        & \Upsilon \frac{1}{m} \sum_{k = 1}^m \| P_k \Delta \theta \|_2 \| P_{\xi} \Delta \theta \|_2 \geqslant \Upsilon \|S \Delta \theta \|_2 \|P_{\xi} \Delta \theta \|_2.
    \end{align*}
\end{proof}
Unfortunately using batching we can not provide gradient estimations satisfying condition~\ref{growth condition}, because $P_k \Delta \theta$ can equals $0$ for some $k$.

\section{Numerical experiments} \label{section exp}
In this section we provide several experiments with an example of real deep learning problems with parameter estimation from definitions~\ref{abs_inexact}, \ref{relative_inexact}, \ref{correl}. We will use Adam optimization algorithm~\cite{kingma2014adam}.
\subsection{Computer Vision} \label{cv exp}
Firstly we examine ResNet18~\cite{he2016deep}:
\begin{figure}[H] \label{resnet18_16}
	\center{\includegraphics[width=1\linewidth]{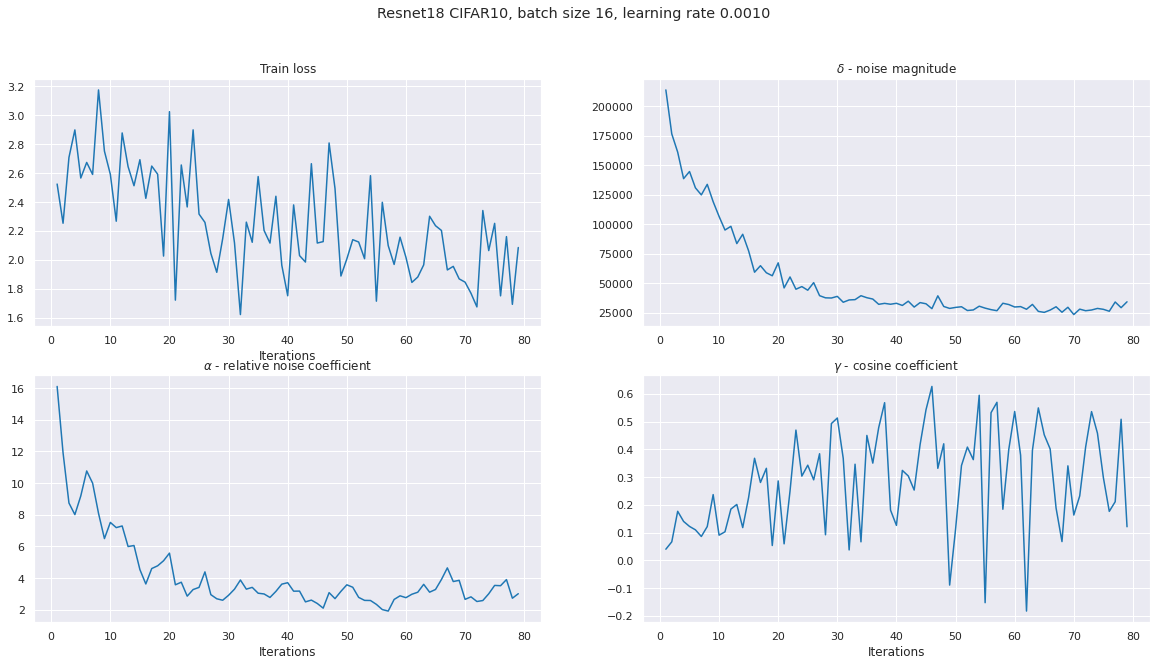}}
	\caption{ResNet-18 CIFAR-10 batch size - 16}
\end{figure}
\begin{figure}[H] \label{resnet18_128}
	\center{\includegraphics[width=1\linewidth]{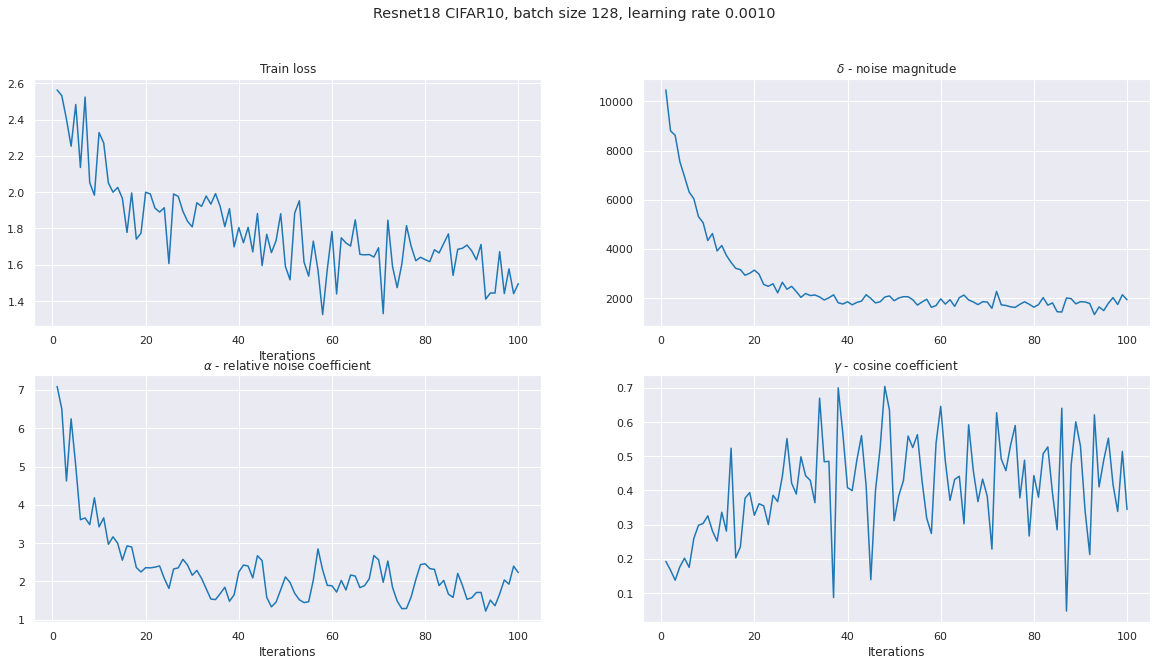}}
	\caption{ResNet-18 CIFAR-10 batch size - 128}
\end{figure}
\begin{figure}[H] \label{resnet18_1024}
	\center{\includegraphics[width=1\linewidth]{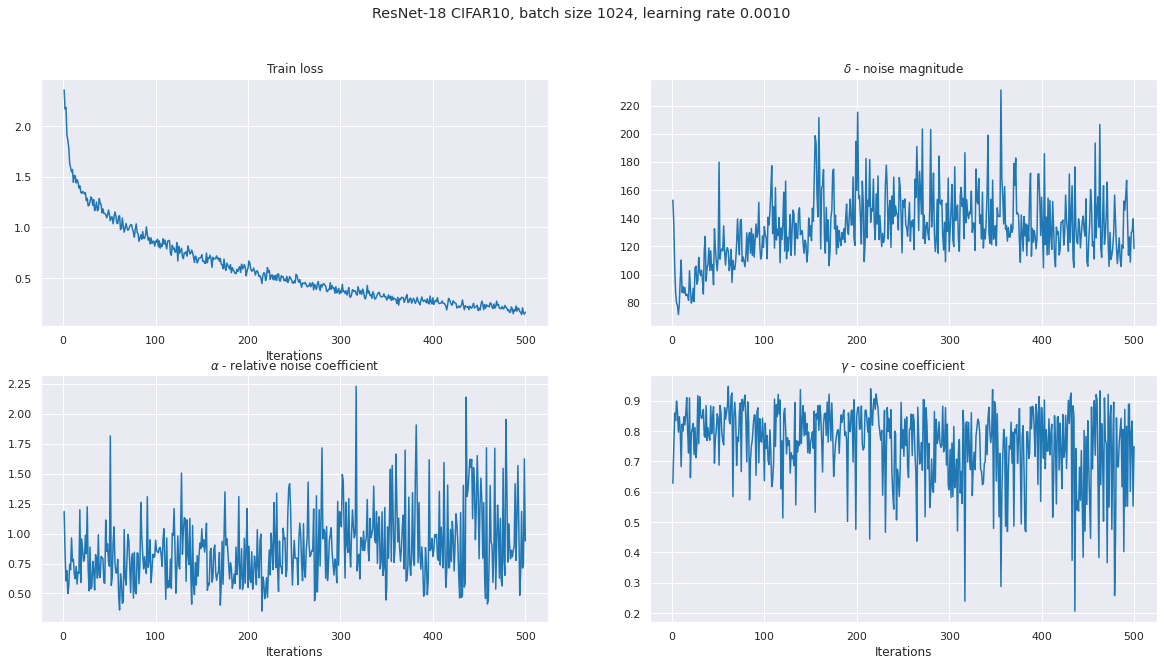}}
	\caption{ResNet-18 CIFAR-10 batch size - 1024}
\end{figure}
%In order to evaluate the impact of residual blocks, consider %InceptionV3~\cite{szegedy2016rethinking}.
%\begin{figure}[H] \label{inception_16}%
%	\center{\includegraphics[width=1\linewidth]{pic_resnet/inception_16.png}}
%	\caption{InceptionV3 CIFAR-10 batch size - 16}
%\end{figure}

\subsection{Natural Language Processing} \label{nlp exp}
Let us examine NLP classification problem on AG News dataset~\cite{zhang2015character}. Firstly consider LSTM~\cite{sak2014long} with 3 layers:
\begin{figure}[H] \label{lstm_1024}
	\center{\includegraphics[width=1\linewidth]{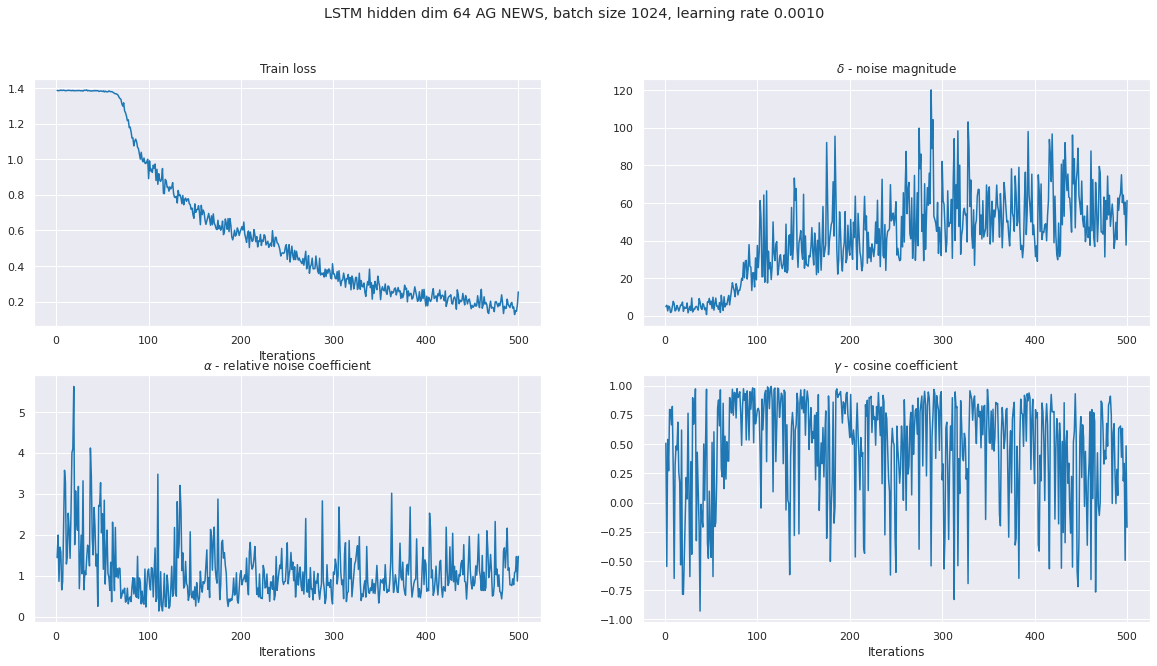}}
	\caption{LSTM AG-News batch size - 1024}
\end{figure}
\begin{figure}[H] \label{lstm_8192}
	\center{\includegraphics[width=1\linewidth]{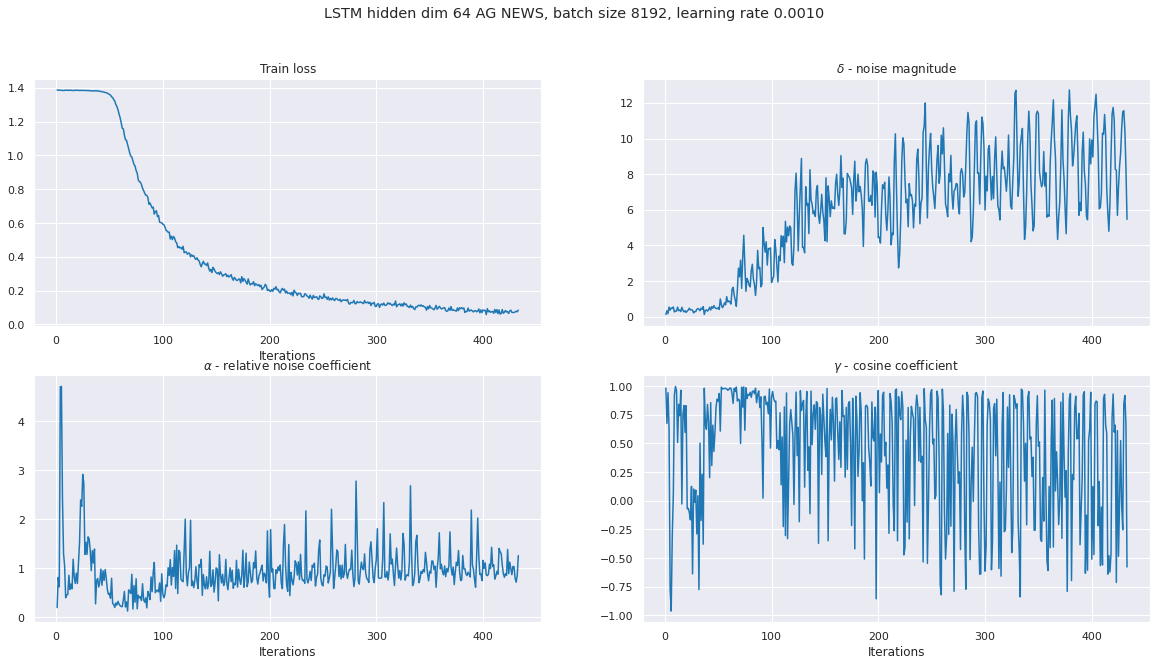}}
	\caption{LSTM AG-News batch size - 8192}
\end{figure}
Then we investigate GRU~\cite{chung2014empirical} with 1 layer:
\begin{figure}[H] \label{gru_128}
	\center{\includegraphics[width=1\linewidth]{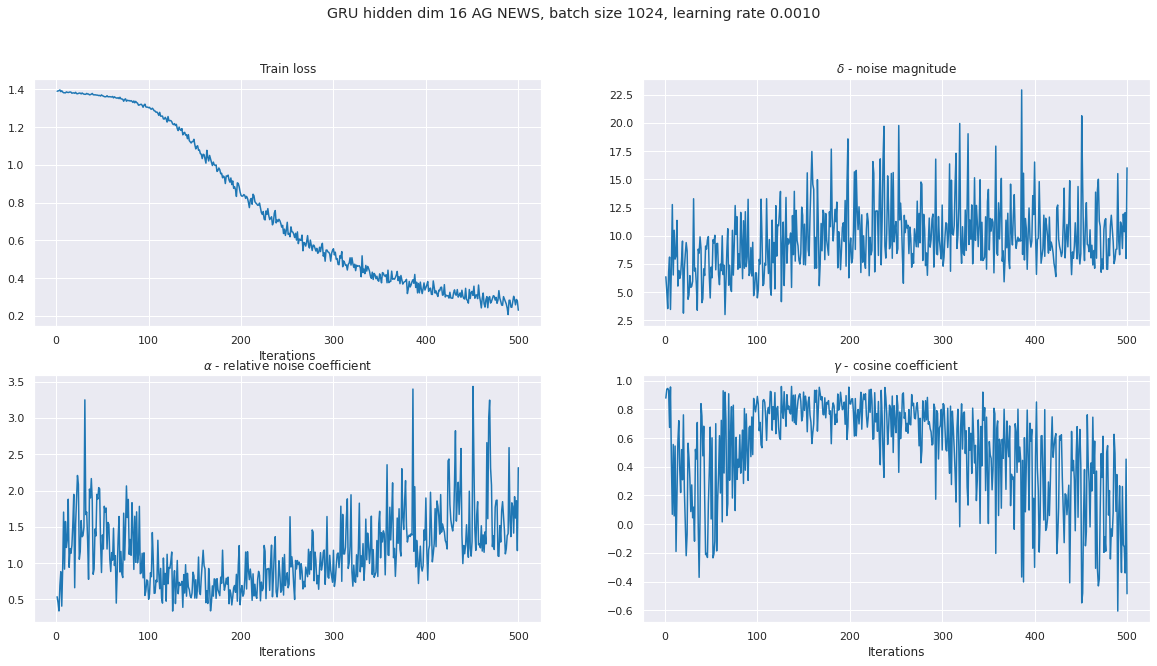}}
	\caption{GRU AG-News batch size - 128}
\end{figure}
\begin{figure}[H] \label{gru_8192}
	\center{\includegraphics[width=1\linewidth]{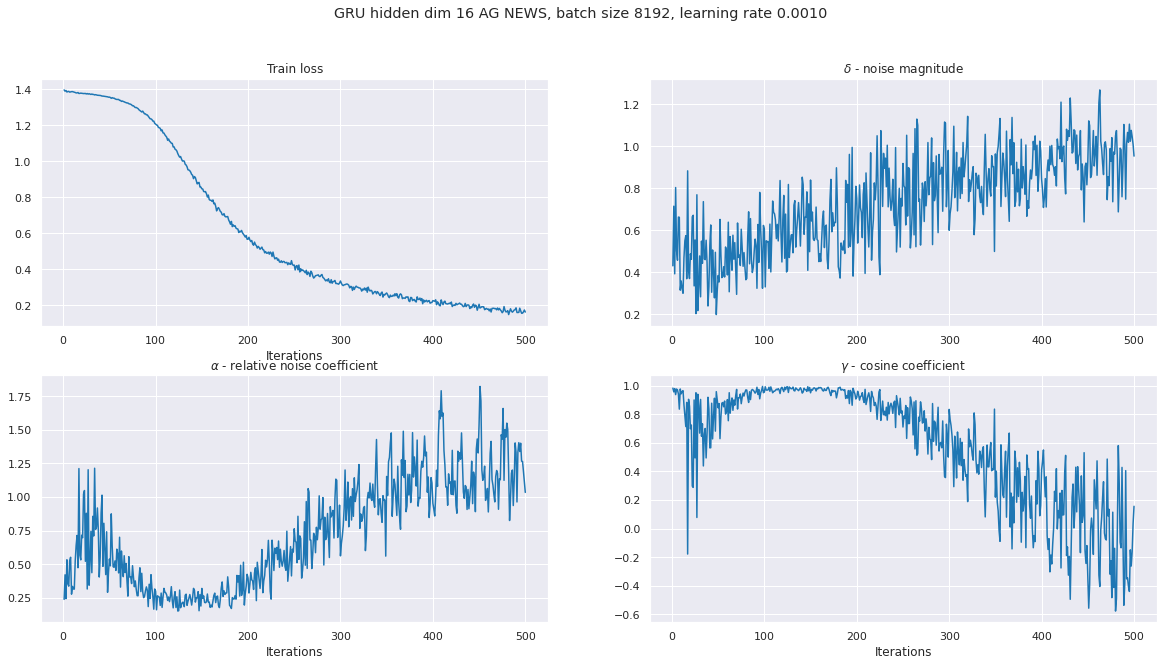}}
	\caption{GRU AG-News batch size - 8192}
\end{figure}

\subsection{Dataset consistency}

We can reformulate condition~\ref{gamma cons} for arbitrary machine learning problem~\ref{sum opt}:
\begin{equation*}
    \begin{gathered}
        f(\theta) = \frac{1}{m} \sum_{k = 1}^m f_k(\theta), \\
        \langle \nabla{f_k}(\theta), \nabla{f_j}(\theta) \rangle \geqslant \Upsilon \|\nabla{f_k}(\theta) \|_2 \| \nabla{f_j}(\theta) \|_2.
    \end{gathered}
\end{equation*}
Such condition may be motivated by consistency of dataset similarly to~\ref{consist data}. Then we can guarantee $\langle \nabla{f}(\theta), \nabla f_{\xi}(\theta) \rangle \geqslant \Upsilon \| \nabla{f}(\theta) \|_2 \| \nabla{f_{\xi}}(\theta) \|_2$. We calculate empirical distribution of 
cosine between two gradient samples on different epochs. We choose random 50 samples and calculate cosine if $i\not=j$.
\begin{figure}[H] \label{cosine_distr}
	\center{\includegraphics[width=1\linewidth]{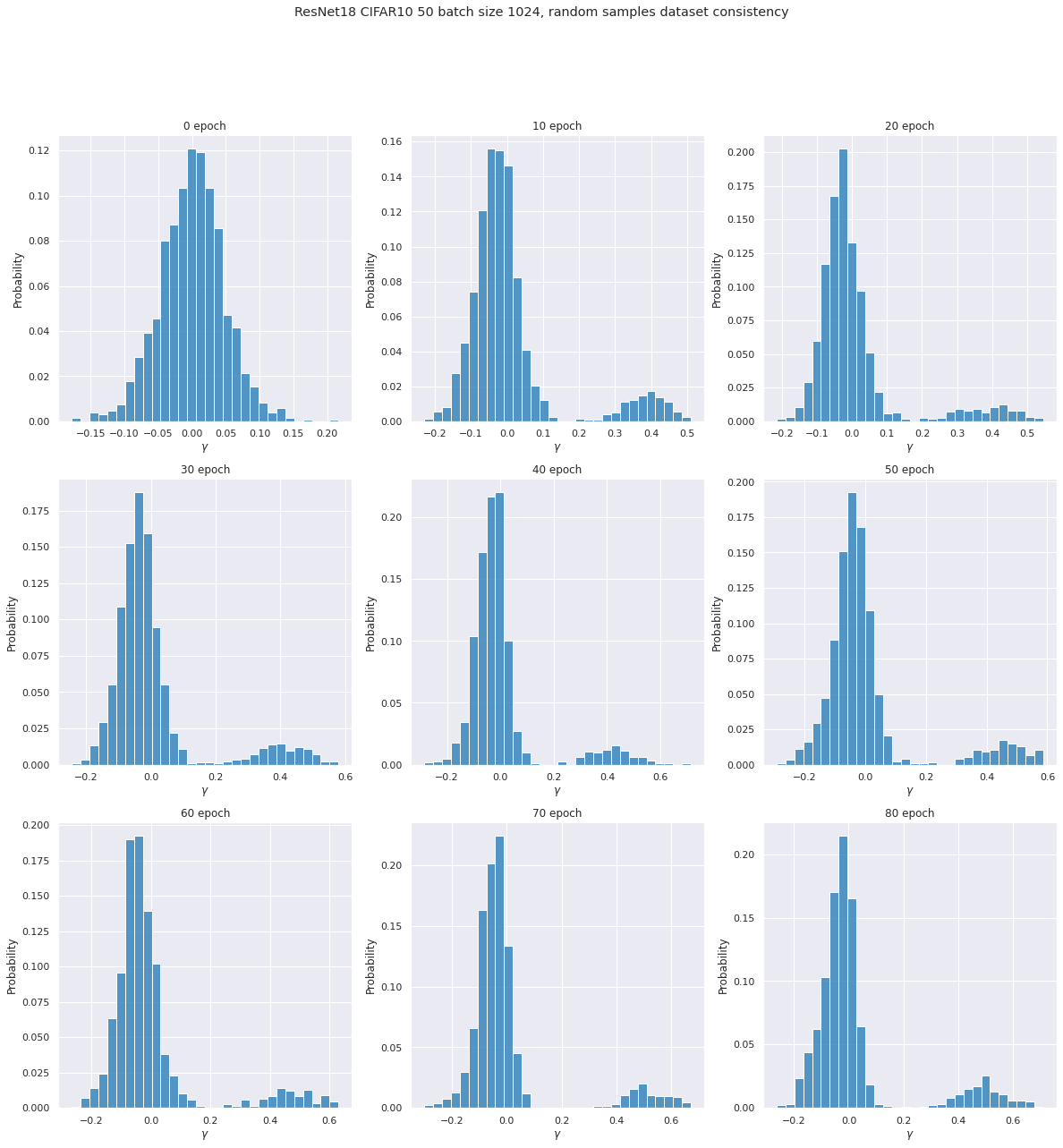}}
	\caption{ResNet-18 CIFAR-10 dataset consistency}
\end{figure}
We can note, that empirical density can be described as mixture of two distributions. First is similar centered normal distribution, and the second one, which can corresponds to group of elements corresponds local convergence of method.

\bibliographystyle{tfs}
\bibliography{literature,all_refs3}
\appendix

\end{document}